\documentclass[10pt,a4paper]{article}
\usepackage[cp1250]{inputenc}
\usepackage[OT1]{fontenc}
\usepackage{amsmath,amssymb}
\usepackage{amsthm}
\usepackage[pdftex]{graphicx}
\usepackage[all]{xy}

\newtheorem{thm}{Theorem}
\newtheorem{lemma}{Lemma}

\newtheorem{conj}{Conjecture}

\newcommand{\laso}{\mathfrak{so}}

\newcommand{\lasl}{\mathfrak{sl}}

\newcommand{\mM}{\textrm M}
\newcommand{\mV}{\textrm V}

\newcommand{\mW}{\textrm W}

\newcommand{\mE}{\textrm E}
\newcommand{\mF}{\textrm F}

\newcommand{\SO}{\textrm{SO}}

\newcommand{\SL}{\textrm{SL}}

\newcommand{\GL}{\textrm{GL}}

\newcommand{\Spin}{\textrm{Spin}}
\newcommand{\G}{\textrm G}

\newcommand{\R}{\mathbb R}

\newcommand{\C}{\mathbb C}
\newcommand{\Sp}{\mathbb S}
\newcommand{\T}{\mathbb T}

\newcommand{\Z}{\mathbb Z}

\newcommand{\cI}{\mathcal I}

\newcommand{\cM}{\mathcal M}

\newcommand{\osu}{\mathcal U}
\newcommand{\osv}{\mathcal V}

\newcommand{\ra}{\rightarrow}

\newcommand{\half}{\frac{1}{2}}

\title{$k$-Dirac operator and Cartan-K\"ahler theorem}
\author{Tom\'a\v s Sala\v c\footnote{Research supported by FWF-Project P23244-N13.}}

\begin{document}
\maketitle
\author
\textit{Abstract: We apply the Cartan-K\"ahler theorem for the k-Dirac operator studied in Clifford analysis and to the parabolic version of this operator. We show that for $k=2$ the tableaux of the first prolongations of these two operators are involutive. This gives us a new characterization of the set of initial conditions for the 2-Dirac operator.}
\section{$k$-Dirac operator}
Let $g$ be the Euclidean product on $\R^{n}$ and let $\{\varepsilon_1,\ldots,\varepsilon_n\}$ be an orthonormal basis. Let $\R_n$ be the Clifford algebra for $(\R^n,g)$ with the defining relation $\varepsilon_\alpha\varepsilon_\beta+\varepsilon_\beta\varepsilon_\alpha=-2g_{\alpha\beta}$. Let $M(n,k,\R)$ be the affine space of matrices of size $n\times k$. Let $\psi$ be a smooth $\R_n$-valued function on $M(n,k,\R)$. We assume throughout this paper that $k\ge2$ and $n\ge3$. For $i=1,\ldots,k$ set
\begin{equation}\label{k-Dirac operator}
\partial_i\psi=\sum_{\alpha=1}^n\varepsilon_\alpha.\partial_{\alpha i}\psi.
\end{equation}
Here $\partial_{\alpha i}$ are the coordinate vector fields on $M(n,k,\R)$. We call the operator $\partial=(\partial_1,\ldots,\partial_k)$ the $k$-\textit{Dirac operator (in the Euclidean setting)} or just the $k$-\textit{Dirac operator}. The $k$-Dirac operator is an overdetermined, constant coefficient system of PDEs. A solution of $\partial\psi=0$ is called a \textit{monogenic function (in the Euclidean setting)} or a \textit{monogenic spinor (in the Euclidean setting)}. More to this operator can be found for example in \cite{CSSS} and \cite{SSSL}. %Stable range means that $n\ge 2k$, otherwise we speak about unstable range.

In this paper we will show that the tableau associated to the first prolongation of the 2-Dirac operator is involutive. This is Theorem \ref{thm involutivity Dirac op}. This gives new characterization of the set of initial conditions. See Theorem \ref{thm initial conditions}. 

We will use representation theory of symmetry group $\SL(k,\R)\times\Spin(n)$ of the operator $\partial$, see \cite{So}. We denote this group by $\G^{ss}_0$. The notation for the symmetry group will be explained afterwards. The group $\G_0^{ss}$ is a semi-simple Lie group with Lie algebra $\lasl(k,\R)\oplus\laso(n)$. We will work with complex representations of $\G^{ss}_0$ and its Lie algebra. We will use the complex spinor representations of $\laso(n)$ rather then the real Clifford module $\R_n$. For $n$ odd, there is only one spinor module $\Sp$. If $n$ is even there are two non-isomorphic spinor modules $\Sp_+$ and $\Sp_-$. In this case we set $\Sp:=\Sp_+\oplus\Sp_-$. We denote by $s$ the dimension of the corresponding spinor module $\Sp$. Thus $s=2^{m}$ if $n=2m+1$ or $n=2m$. 

The question of finding the set of initial conditions for the $k$-Dirac operator might be the right link which brings into the game a parabolic $k$-Dirac operator. This link is explained in the next paragraph. 

\subsection{Parabolic $k$-Dirac operator}\label{section parabolic D operator}
The \textit{parabolic k-Dirac operator} $D$ is an invariant first order operator which lives in the world of parabolic geometries. For the purpose of this paper we will give a coordinate definition of this operator. Put $\osu=M(n,k,\R)\times A(k,\R)$ where $A(k,\R)$ is the affine space of skew-symmetric matrices of size $k$. We write coordinates as $(x_{\alpha i},y_{rs})$ where $x_{\alpha i}$, resp. $y_{rs}$ are coordinates on $M(n,k,\R)$, resp. on $A(k,\R)$. We write $\partial_{\alpha i}=\partial_{x_{\alpha i}},\partial_{rs}=\partial_{y_{rs}}$. We use the convention that $\partial_{rs}=-\partial_{sr}$. The set $\osu$ is a isomorphic to an open affine subset of the Grassmannian of isotropic $k$-planes in $\R^{k,n+k}$. The Grassmannian is a flat model for a particular type of parabolic geometries. For $k=2$ this geometry is refered to as Lie contact structures, see \cite{CSl}.

For $\alpha=1,\ldots,n$ and $i=1,\ldots,k$ put $L_{\alpha i}:=\partial_{\alpha i}-\frac{1}{2}x_{\alpha j}\partial_{ji}$. We will call these vector fields \textit{left invariant vector fields}. Lie bracket is 
\begin{eqnarray}\label{Lie bracket of left invariant vector fields}
[L_{\alpha i},L_{\beta j}]=g_{\alpha\beta}\partial_{ij}.
\end{eqnarray}
These vector fields span a non-integrable distribution on $\osu$ which is the essence of the parabolic geometry. The (graded) tangent bundle of the Grassmannian variety has a natural reduction of the structure group to $\GL(k,\R)\times\SO(n)$. We may lift (uniquely) the trivial principal $\GL(k,\R)\times\SO(n)$-bundle over $\osu$ to $\G_0=\GL(k,\R)\times\Spin(n)$-bundle. This is the usual notation from \cite{CSl}. The group $\G_0$ is a reductive group whose semi-simple part is isomorphic to $\G_0^{ss}$. We extend the action of $\G_0^{ss}$ on $\Sp$ to the action of $\G_0$ by the choice of a generalised conformal weight, i.e. we specify the action of the center of $\GL(k,\R)$, as in \cite{St}. By associating $\Sp$ to the principal $\G_0$-bundle we get a spinor bundle over $\osu$.

We may view the set of vector fields $\{L_{\alpha i},\partial_{rs}|\alpha=1,\ldots,n;i=1,\ldots,k;1\le r<s\le k\}$ as a section of the principal $\GL(k,\R)\times\SO(n)$-bundle. This will be over preferred gauge over $\osu$. We choose a section of the $\G_0$-bundle compatible with the preferred gauge. Then a section of the spinor bundle $\psi$ becomes a spinor valued function on $\osu$. Then with the choice of the very flat Weyl connection we can write in this gauge $D\psi=(D_1\psi,\ldots,D_k\psi)$ where
\begin{equation}\label{parabolic k-Dirac operator}
D_i\psi=\sum_{\alpha=1}^n\varepsilon_\alpha. L_{\alpha i}\psi.
\end{equation}
Comparing this to (\ref{k-Dirac operator}) we see that we have just replaced each $\partial_{\alpha i}$ by the corresponding left invariant vector field $L_{\alpha i}$. A solution of $D\psi=0$ is called a \textit{(parabolic) monogenic spinor}. 

There is a strong and very interesting link which leads from the operator (\ref{parabolic k-Dirac operator}) to (\ref{k-Dirac operator}). First of all, a parabolic monogenic spinor $\psi$ which does not depend on $y$-coordinates can be naturally viewed as a solution of $\partial\psi=0$. A bit of work shows that there is a (unique) locally exact sequence of invariant operators starting with the operator $D$. We can do the same move for the whole sequence as we did with the monogenic spinors, i.e. we work only with the real analytic sections which in the preferred trivialization do not depend on $y$-coordinates on which we act by the invariant operators. Then we get a new sequence of operators which is still locally exact and thus descends to a resolution of $\partial$. This can be found in  \cite{S}.

Motivation for this paper is hidden in the set of initial conditions for these two systems of PDEs. It is not hard to see that any monogenic spinor (in the Euclidean setting) $\psi$ is uniquely determined by its restriction to the set $M(n-1,k,\R)\cong\{x_{11}=x_{12}=\ldots=x_{1k}=0\}$. Moreover on this set the restriction $\psi|_{M(n-1,k,\R)}$ has to satisfy for each $i,j=1,\ldots,k:$ 
\begin{equation}\label{second order equation}
[\tilde{\partial_i},\tilde{\partial_j}]\psi|_{M(n-1,k,\R)}=0
\end{equation}
where $\tilde{\partial}_i=\sum_{\alpha=2}^n\varepsilon_\alpha\partial_{\alpha i}$. This is a consequence of the fact that the coordinate vector fields commute. On the other hand given a real analytic spinor valued function $\varphi$ on $M(n-1,k,\R)$ converging on open neighbouhood of $x\in M(n-1,k,\R)$ and which satisfies (\ref{second order equation}) then there is a unique monogenic spinor on $M(n,k,\R)$ converging on some open neighbouhood of $x$ whose restriction to $M(n-1,k,\R)$ coincides with $\varphi$. 

\begin{conj}
Given arbitrary real analytic spinor $\psi$ in $x_{\alpha i}$-variables with $\alpha\ge2$ converging on some open subset $\osv$ of $M(n-1,k,\R)$ there is a unique (parabolic) monogenic spinor $\Psi$, i.e. $D\Psi=0$, convering on some open neighbouhood of $\osv$ in $\osu$ whose restriction to the set $M(n-1,k,\R)\cong\{x_{11}=\ldots=x_{1k}=y_{12}\ldots=y_{k-1,k}=0\}$ coincides $\psi$.
\end{conj}

If Conjecture 1. is true then the system (\ref{parabolic k-Dirac operator}) will have a nicer set of initial conditions then (\ref{k-Dirac operator}). Starting with the $k$-Dirac operator in the Euclidean setting and looking for a new system of PDEs such that any quadratic real analytic spinor on $M(n-1,k,\R)$ extends to a unique solution of the new system then one can derive the Lie bracket (\ref{Lie bracket of left invariant vector fields}) and the right dimension of the set $\osu$. This is already a link from the operator (\ref{k-Dirac operator}) to the operator (\ref{parabolic k-Dirac operator}). The only question is how good this link is?

I hoped this result would follow from the Cartan-K\"ahler theorem. Unfortunatelly it does not. Nevertheless the Cartan-K\"ahler theorem gives us some other interesting results. In this paper we show that both systems are involutive after the first prolongation if $k=2$. These are Theorems \ref{thm involutivity Dirac op}. and \ref{thm involutivity D op}. For $k\ge 3$ this is no longer true and one has to continue on prolongating. I do not know when the involutivity is attained. A closer look on the proof of involutivity for the parabolic 2-Dirac operator also explains why the Cartan-K\"ahler theorem does not give Conjecture 1.

In the next section we cover basic machinery and terminology needed for the Cartan-K\"ahler theorem. This short summary is taken mostly from \cite{IL}. For more on the Cartan-K\"ahler theorem and exterior differential systems see \cite{B}. 

If it is from the context clear whether we talk about $D$ or $\partial$ then we will simply say a $k$-Dirac operator. Similarly we say just a monogenic spinor $\psi$ if it clear if we mean $\partial\psi=0$ or $D\psi=0$.

\section{Exterior differential systems}
We assume all over the paper that all structures are real analytic. Then we can apply machinery of the Cartan-K\"ahler theorem.

Let $M$ be a manifold. An \textit{exterior differential system} (EDS) on $M$ is a graded differential ideal $\cI\subset\Omega^\ast(M)$. Recall that $\Omega^\ast(M)$ is naturally graded by the degree of differential forms. A graded differential ideal is a graded ideal closed under the de Rham differential. We denote the $k$-th homogeneous piece by $\cI^k$. We are interested in \textit{integral manifolds} of $\cI$. These are submanifolds $i:N\hookrightarrow M$ such that $i^\ast\alpha=0$ for any $\alpha\in\cI$. Many interesting problems can be formulated in the language of EDSs and integral manifolds. For a fixed $x\in M$, the set of \textit{integral elements} for $\cI$ of rank $k$ at $x$ is the set $\{E\subset T_xM:\dim(E)=k,\forall\alpha\in\cI^k:\alpha|_E=0\}$ 

We will be interested in EDSs with independence condition. An independence condition for $\cI$ is given by a set of 1-forms $J=\{\omega^1,\ldots,\omega^n\}$. We consider only those integral manifolds of $\cI$ for which $i^\ast(\omega^1\wedge\ldots\wedge\omega^n)$ is a non-vanishing form on $N$. 

Let $I$ be the ideal generated by forms $\{\theta^1,\ldots,\theta^s\}$ such that these forms generate $\cI$ as a differential ideal. We will be interested here only in the case when each $\theta^i$ is a 1-form. We call such system a \textit{Pfaffian system}. We may assume that $\{\omega^1,\ldots,\omega^n,\theta^1,\ldots,\theta^s\}$ is a set of everywhere linearly independent 1-forms on $M$. Let $\pi^1,\ldots,\pi^t$ be a set of 1-forms such that the set $\{\omega^i,\theta^j,\pi^\varepsilon\}$ with $i=1,\ldots,n,j=1,\ldots,s,\varepsilon=1\ldots,r$ is a basis of $T^\ast_x M$ for each $x\in M$. Let us now fix $x\in M$. We set $\mV^\ast:=(J/I)_x,\mW^\ast:=I_x$. We write $v^i=\omega^i_x,w^j=\theta^j_x$ and denote the dual elements by $v_i$ and $w_j$.

Then the Pfaffian system $\cI$ is called a \textit{linear Pfaffian system} if 
\begin{equation}\label{structure equation1}
d\theta^a=A^a_{\varepsilon i}\pi^\varepsilon\wedge\omega^i+T^a_{ij}\omega^i\wedge\omega^j\ mod\ I
\end{equation}
holds for some functions $A^a_{\varepsilon b},T^a_{ij}$. The \textit{tableau} at the point $x$ is equal to  $A_x:=\{A^a_{\varepsilon i}v^i\otimes w_a\subset\mV^\ast\otimes\mW|1\le\varepsilon\le r\}$. We drop the subscript $x$ and write $A$ instead of $A_x$. 

Let $\delta:\mV^\ast\otimes\mV^\ast\otimes\mW\ra\Lambda^2\mV^\ast\otimes\mW$ be the natural projection. Set $H^{0,2}(A):=\Lambda^2\mV^\ast\otimes\mW/\delta(\mV^\ast\otimes A)$. The \textit{torsion} of $(I,J)$ at $x$ is defined as the class $[T]_x:=[T^a_{ij}(x) v^i\wedge v^j\otimes w_a]\in H^{0,2}(A)$. If $[T_x]=0$ we say that the torsion is \textit{absorable} at x. This means that we can replace $\pi^\varepsilon$ by new forms $\pi'^\epsilon,\epsilon=1,\ldots,t$ such that $\{\omega^i,\theta^j,\pi'^\varepsilon\}$ is still a basis of $T^\ast_xM$ for each $x\in M$ and $T'^a_{ij}=0$ for all $a,i,j$ at the point $x$ with respect to the new basis. Later on we will need the following lemma.

\begin{lemma}\label{lemma vanishing of torsion}
For a linear Pfaffian system on $M$ with an independence condition the necessary and sufficient condition for vanishing of the torsion $[T_x]$ at a point $x\in M$ is that there is an integral element over $x\in M$ satisfying the independence condition.
\end{lemma}
Proof: See Proposition 5.14. from \cite{B}.

We see that the torsion is a primary obstruction for existence of integral manifolds. Suppose that $[T]=0$ on an open neighbouhood of $x$. Choose a basis $\{u^1,\ldots,u^n\}$ of $\mV^\ast$ and set for $k=1,\ldots,n: A_k:=A\cap\langle u^{k+1},\ldots,u^n\rangle\otimes\mW$. Here $\langle\ \rangle$ denotes the linear span. Put $A^{(1)}:=S^2\mV^\ast\otimes\mW\cap\mV^\ast\otimes A$. We call $A^{(1)}$ the \textit{(first) prolongation of the tableau} $A$. Then we have the inequality
\begin{equation}\label{Cartan test}
\dim(A^{(1)})\le\dim(A)+\dim(A_1)+\ldots+\dim(A_{n-1}).
\end{equation} 
We say that the tableau is \textit{involutive} if the equality holds for some choice of a basis of $\mV^\ast$. It is convenient to introduce the \textit{Cartan characters} $s_1,\ldots,s_n$ of the tableau by requiring that $\dim(A)-\dim(A_k)=s_1+\ldots+s_k$ holds for each $k=1,\ldots,n$. Then the inequality (\ref{Cartan test}) becomes
\begin{equation}\label{Cartan test1}
\dim(A^{(1)})\le s_1+2s_2+\ldots+ns_n.
\end{equation}

If the tableau is involutive then the Cartan-K\"ahler theorem applies. The Cartan-K\"ahler theorem guarantees existence of $n$-dimensional integral manifolds passing through $x$ satisfying the independence condition. Moreover we can read from the Cartan characters "how many" such local manifolds there are.

If the tableau is not involutive one has to start over on the pullback of the canonical system on the Grassmann bundle to the space of integral elements. In calculation this means that we add elements from $A^{(1)}$ as new variables and add new forms $\theta^a_{i}:=A^a_{\varepsilon i}\pi^\varepsilon-p^a_{ij}\omega^j$ where $p^a_{ij}v^i\otimes v^j\otimes w_a\in A^{(1)}$ to the ideal $I$.

EDSs naturally arise with PDEs. Suppose that we are given a system of PDEs of order $k$. Then we take $M$ to be the space of $k$-jets of solutions of the PDE and we pull back the canonical system which lives on the space of jets of vector valued functions. 

Computation simplifies in the case of a constant coefficient system. The torsion vanishes and in the case of a linear, constant coefficient, homogeneous system of PDEs the tableau $A$ is at any point isomorphic to the space of linear solutions of this system. The first prolongation of the tableau is naturally isomorphic to the space of quadratic solutions. Set $A^{(0)}:=A$ and inductively for $j=1,2,\ldots$ put $A^{(j)}:=S^j\mV^\ast\otimes\mW\cap\mV^\ast\otimes A^{(j-1)}$. Then $A^{(j)}$ is naturally isomorphic the to space of homogeneous solutions of the system of homogeneity $j+1$ and $A^{(j+1)}\cong (A^{(j)})^{(1)}$. For more see \cite{IL}.

\section{$k$-Dirac operator (in the Euclidean setting) and the Cartan-K\"ahler theorem}\label{section CK for k-Dirac op}

%\subsection{Homogeneous monogenic spinors}

For this paper we will need to understand the space of linear, quadratic and cubic monogenic spinors. Recall that $\Sp$ is the complex spinor representation of $\laso(n)$ defined in Section 1. There is an isomorphism $M(k,n,\R)\otimes_\R\Sp\cong M(k,n,\C)\otimes_\C\Sp$. We will work with complex representations of Lie algebra of $\G_0^{ss}$ and take tensor product over complex numbers. We will denote the Cartan product by $\boxtimes$. This is the irreducible subspace with the highest weight in the tensor product of irreducible representations.

Let $\mE$, resp. $\mF$ be the defining representation of $\lasl(k,\C)$, resp. of $\laso(n,\C)$. We choose a basis $\{e_1,\ldots,e_k\}$ of $\mE$. We denote by $\{\varepsilon_1,\ldots,\varepsilon_n\}$ an orhonormal basis of $\mF$. Let $g$ be the $\SO(n,\C)$-invariant scalar product on $\mF$. If $n=2m$ is even we denote by $\{v_1,\ldots,v_m,w_1,\ldots,w_m\}$ a null basis of $\mF$ such that $g(v_i,w_j)=\delta_{ij}$. If $n=2m+1$ is odd we denote by $\{v_1,\ldots,v_m,w_1,\ldots,w_m,u\}$ a basis such that the relations on $v_i,w_j$ are same as for $n=2m$ and $g(u,u)=1,g(u,w_i)=g(u,v_j)=0$.

In the case of the $k$-Dirac operator we have that $(\mV_x)_\C:=(\mV^\ast_x)\otimes_\R\C\cong\mE\otimes\mF$ and $\mW_x\cong\Sp$ for each $x\in M(n,k,\R)$. The tableau is isomorphic to $A_x\cong\mE\otimes\T$ where $\T\subset\mF\otimes\Sp$ is the Cartan component. We will drop the subscript $x$. With the notation above $\T\cong\mF\boxtimes\Sp$. We call $\T$ the \textit{twistor representation} of $\laso(n)$. The subspace $\T$ is invariantly defined as the kernel of the canonical projection $\pi:\mF\otimes\Sp\ra\Sp$ which is on decomposable elements given by the Clifford multiplication $\pi(\varepsilon\otimes s)=\varepsilon .s$. By induction on $i$ we get that $A^{(i)}$ is the intersection of $S^{i+1}(\mE\otimes\mF)\otimes\Sp$ with the kernel of the projection 
\begin{equation}\label{description of prolongation}
\mE^{\otimes^{i+1}}\otimes\mF^{\otimes^{i+1}}\otimes\Sp\xrightarrow{Id_\mM\otimes\pi}\mE^{\otimes^{i+1}}\otimes\mF^{\otimes^{i}}\otimes\Sp
\end{equation}
where $\mM=\mE^{\otimes^{i+1}}\otimes\mF^{\otimes^i}$.

\subsection{The space of polynomials on $M(n,k,\C)$ as a $\GL(k,\C)\times\GL(n,\C)$-module}\label{section on decomposition of polynomials}
Let us consider the action of $\GL(k,\C)\times\GL(n,\C)$ on the space of polynomials on $M(n,k,\C)$ given by $((g,h).f)(x)=f(hxg^T)$ where $g\in\GL(k,\C),h\in\GL(n,\C),x\in M(n,k,\C)$ and $f$ is a polynomial on $M(n,k,\C)$. The space of linear polynomials is isomorphic to $\mE'\otimes\mF'$ where $\mE',\mF'$ is the defining representation of $\GL(k,\C)$, resp. of $\GL(n,\C)$. The space of quadratic polynomials is then isomorphic to $S^2(\mE'\otimes\mF')\cong S^2\mE'\otimes S^2\mF'\oplus\Lambda^2\mE'\otimes\Lambda^2\mF'$ The set $\{e_i\odot e_j\otimes\varepsilon_\alpha\odot\varepsilon_\beta,e_i\wedge e_j\otimes\varepsilon_\alpha\wedge\varepsilon_\beta|i,j=1,\ldots,k;\alpha,\beta=1,\ldots,n\}$ is a basis of $S^2(\mE'\otimes\mF')$. Here we are using the bases introduced in the previous paragraph. The corresponding polynomials are $x_{\alpha i}x_{\beta j}+x_{\beta i}x_{\alpha j}\in S^2\mE'\otimes S^2\mF'$, resp. $x_{\alpha i}x_{\beta j}-x_{\beta i}x_{\alpha j}\in\Lambda^2\mE'\otimes\Lambda^2\mF'$. 

With respect to the usual choices on Lie algebra of the semi-simple part of $\GL(k,\C)\times\GL(n,\C)$, i.e. the Cartan subalgebra consists of diagonal matrices and positive roots span the strictly upper triangular matrices,
the polynomials $x_{11},x_{11}x_{22}-x_{12}x_{21}$ are highest weight vectors of $\mE'\otimes\mF'$, resp. of $\Lambda^2\mE'\otimes\Lambda^2\mF'$. If $k=2$ then Theorem 5.2.7 on $\GL(k,\C)\times\GL(n,\C)$-duality from \cite{GW} shows that any highest weight polynomial of an $\GL(k,\C)\times\GL(n,\C)$-irreducible subspace is up to a scalar multiple a product of $x_{11}$ and $x_{11}x_{22}-x_{12}x_{21}$. 

\subsection{Non-involutivity of the tableau of the $k$-Dirac operator} 
Recall that we have denoted by $\{\varepsilon_1,\ldots,\varepsilon_n\}$ an orthonormal basis of $\C^n$. We consider it also as an orthonormal basis of $\R^n$ and $\{e_1,\ldots,e_k\}$ as a basis of $\R^k$. Then $\{e_i\otimes\varepsilon_\alpha\}$ is a basis of $\mV^\ast$. Let us order it by $\{e_1\otimes\varepsilon_1,\ldots,e_1\otimes\varepsilon_k,e_2\otimes\varepsilon_1,\ldots,e_n\otimes\varepsilon_k\}$. Let $\{s^\mu|\mu=1,\ldots,s\}$ be a basis of $\Sp$. Then $\{\varepsilon_\alpha\otimes s^\mu+\varepsilon_n\otimes\varepsilon_n.\varepsilon_\alpha.s^\mu|\alpha<n\}$ is a basis of $\T$. In particular $\dim(\T)=s(n-1)$. The Cartan characters with respect to the ordered basis of $\mV^\ast$ are equal to $s_1=\ldots=s_{k(n-1)}=s,s_{k(n-1)+1}=\ldots=s_{nk}=0$. It is clear that these characters minimaze the right hand side of (\ref{Cartan test1}). So the right hand side in (\ref{Cartan test1}) is equal to $s{k(n-1)+1\choose2}$. 

We now have to compute the dimension $\dim(A^{(1)})$. We know that $A^{(1)}$ is naturally isomorphic to the space of quadratic monogenic spinors. %Let $\{s^\mu|\mu=1,\ldots,s\}$ be a basis of $\Sp$.

\begin{lemma}\label{quadratic monogenic spinors lemma}
The dimension of $A^{(1)}$ is equal to $s{k(n-1)+1\choose 2}-s{k\choose 2}$. In particular $\dim A^{(1)}<s{k(n-1)+1\choose 2}$ and the tableau associated to the $k$-Dirac operator is not involutive.
\end{lemma}
Proof: First let us consider the piece $A^{(1)}\cap S^2\mE\otimes S^2\mF\otimes\Sp$. Let $a^\mu_{ij\alpha\beta}\in\C$ with $\alpha,\beta>1$ be artibrary such that $a^\mu_{ij\alpha\beta}=a^\mu_{ij\beta\alpha}=a^\mu_{ji\alpha\beta}$. Consider the element
%\begin{eqnarray*}
$a_{ij\alpha\beta}^\mu e_i\odot e_j\otimes(\varepsilon_\alpha\odot\varepsilon_\beta\otimes s^\mu+\varepsilon_\alpha\odot\varepsilon_1\otimes\varepsilon_1.\varepsilon_\beta. s^\mu
+\varepsilon_\beta\odot\varepsilon_1\otimes\varepsilon_1.\varepsilon_\alpha. s^\mu-2\delta^{\alpha\beta}\varepsilon_1\otimes\varepsilon_1\otimes s^\mu).$
%\end{eqnarray*}
Then from (\ref{description of prolongation}) follows that this is an element of $S^2\mE\otimes S^2\mF\otimes\Sp\cap A^{(1)}$ with the given term $a_{ij\alpha\beta}^\mu e_i\odot e_j\otimes\varepsilon_\alpha\odot\varepsilon_\beta\otimes s^\mu$ where $\alpha,\beta>1$. On the other hand any element of $S^2\mE\otimes S^2\mF\otimes\Sp\cap A^{(1)}$ is uniquely determined by the coefficients $a^\mu_{ij\alpha\beta}\in\C$ with $\alpha,\beta>1$. This shows that $\dim(S^2\mE\otimes S^2\mF\otimes\Sp\cap A^{(1)})=s{n\choose2}{k+1\choose2}$. As a $\lasl(k,\C)\oplus\laso(n,\C)$-module this piece is isomorphic to $S^2\mE\otimes S^2_0\mF\boxtimes\Sp$ where $S^2_0\mF$ denotes the trace-free part of $S^2\mF$.

Let $a^\mu_{ij\alpha\beta}\in\C$ with $\alpha,\beta>1$ be such that $a^\mu_{ij\alpha\beta}=-a^\mu_{ij\beta\alpha}=-a^\mu_{ji\alpha\beta}$. Then the element
%\begin{eqnarray*}
$a_{ij\alpha\beta}^\mu e_i\wedge e_j\otimes(\varepsilon_\alpha\wedge\varepsilon_\beta\otimes s^\mu+\varepsilon_\alpha\wedge\varepsilon_1\otimes\varepsilon_1.\varepsilon_\beta. s^\mu
-\varepsilon_\beta\wedge\varepsilon_1\otimes\varepsilon_1.\varepsilon_\alpha.s^\mu)$
%\end{eqnarray*}
belongs to $A^{(1)}$ iff for all $i,j=1,\ldots,k:\sum_{\alpha\beta}a^\mu_{ij\alpha\beta}\varepsilon_\alpha.\varepsilon_\beta.s^\mu=0$. These are $s{k\choose2}$ linearly independent equations. This shows that $\dim(\Lambda^2\mE\otimes\Lambda^2\mF\otimes\Sp\cap A^{(1)})=s({k\choose 2}{n -1\choose 2}-{k\choose 2})$. As a $\lasl(k,\C)\oplus\laso(n,\C)$-module this piece is isomorphic to $\Lambda^2\mE\otimes\Lambda^2\mF\boxtimes\Sp$. This space is irreducible if $n>4$. If $n=4$ it is the direct sum of two irreducible pieces. Summing up we get that
$\dim(A^{(1)}) =%&=&s{n\choose2}{k+1\choose2}+s({k\choose 2}{n -1\choose 2}-{k\choose 2})\\
%&=&\frac{s}{4}k(n-1)(n(k+1)+(k-1)(n-2))-s{k\choose 2}\\
%&=&\frac{s}{2}k(n-1)(nk-k+1)-s{k\choose 2}\\
s{k(n-1)+1\choose 2}-s{k\choose 2}. \Box$
%\end{eqnarray*}

\subsection{Involutivity of the tableau of the first prolongation for $k=2$}
Now we show that the tableau associated to the first prolongation is involutive when $k=2$. We replace $A$ by $A^{(1)}$ and $A^{(1)}$ by $A^{(2)}$ and repeat the algorithm. This means that we have to compare the sum of the Cartan characters with respect to a suitable filtration on the space of quadratic monogenic spinors (right hand side of (\ref{Cartan test1})) to the dimension of the space of cubic monogenic spinors. 

According to Section \ref{section on decomposition of polynomials} the space of homogeneous polynomials of degree 3 on $M(n,2,\C)$ decomposes into $S^3\mE'\otimes S^3\mF'\oplus\mE'\boxtimes\Lambda^2\mE'\otimes(\mF'\boxtimes\Lambda^2\mF')$ as a $\GL(2,\C)\times\GL(n,\C)$-module. %The corresponding highest weight polynomials are $x_{11}^3,x_{11}(x_{11}x_{22}-x_{12}x_{21})$.
If we restrict $S^3\mF'$ to $\laso(n,\C)$ then it decomposes into $S^3_0\mF\oplus\mF$ where $S^3_0\mF$ is the trace-free part of $S^3\mF$. The trace-free part is the kernel of the canonical contraction. If $n>4$ then $\Lambda^2\mF'$ is irreducible also under the action of $\laso(n,\C)$. If $n=4$ then $\laso(4,\C)\cong\laso(3,\C)\oplus\laso(3,\C)$ and the defining representation $\C^4$ is isomorphic $\C^4\cong\C^2\otimes\C^2$ where we use that $\laso(3,\C)\cong\lasl(2,\C)$. The space $\Lambda^2\C^4$ decomposes into $\Lambda^2\C^2\otimes S^2\C^2\oplus S^2\C^2\otimes\Lambda^2\C^2$. This is a splitting of 2-forms into self-dual and anti-self-dual part. The spinor representations are $\Sp_+\cong\C^2\otimes\C,\Sp_-\cong\C\otimes\C^2$, i.e. it is the defining representation of one summand times the trivial representation of the second summand. The projection $\pi:(\C^2\otimes\C^2)\otimes\Sp_+\ra\Sp_-$ is then the obvious skew-symmetrization in the first factor times the identity on the latter factor.  Similarly for $\Sp_-$.

\begin{lemma}\label{lemma second prolongation}
Let us write $n=2m$ if $n$ is even and $n=2m+1$ if $n$ is odd. The space of cubic monogenic spinors contains the following list of irreducible $\lasl(2,\C)\oplus\laso(n,\C)$-modules
\begin{center}\label{ternary monogenic spinors}
\begin{tabular}{llll}
$n\ge 3:$&$S^3\mE\otimes S^3_0\mF\boxtimes\Sp$&$\frac{2^{m+1}(n+1)n(n-1)}{3}$\\
%$n=4:$& $(2,1)\otimes\mM$&&\\
$n\ge 5:$&$\mE\boxtimes\Lambda^2\mE\otimes\Lambda^2\mF\boxtimes\mF\boxtimes\Sp$&$\frac{2^{m-1}(n+1)(n-1)(n-3)}{3}$
\end{tabular}
\end{center}
where on the right is the dimension of the module. For $n=4$ there is also the module $\mE\boxtimes\Lambda^2\mE\otimes(\C^2\otimes S^4\C^2\oplus S^4\C^2\otimes\C^2)$.
\end{lemma}
Proof: Let $v_1\in\mF$ be the first basis element introduced in Section \ref{section CK for k-Dirac op}. Then $v_1$ is null and we may assume that it is a highest weight vector in $\mF$. Let $s\in\Sp$ be a highest weight vector. Then $v_1.s=0$. This shows that $S^3\mE\otimes S^3\mF\boxtimes\Sp\subset A^{(1)}$ if $n\ge 3$. If $n=4$ then $\mF\otimes\Lambda^2\mF$ %\cong$ (\C^2\otimes\C^2)\otimes(\Lambda^2\C^2\otimes S^2\C^2\oplus S^2\C^2\otimes\Lambda^2\C^2)$.  $(\C^2\otimes\C^2)\boxtimes(\Lambda^2\C^2\otimes S^2\C^2)\cong
contains two unique pieces $\C^2\otimes S^3\C^2\oplus%,(\C^2\otimes\C^2)\boxtimes(S^2\C^2\otimes\Lambda^2\C^2)\cong 
S^3\C^2\otimes\C^2$. Then $S^3\C^2\otimes\C^2\boxtimes\Sp_+\cong S^4\C^2\otimes\C^2$ is in the kernel of the map (\ref{description of prolongation}). Similarly for $\Sp_-$. This shows that $(2,1)\otimes(\C^2\otimes S^4\C^2\oplus S^4\C^2\otimes\C^2)$ belongs to $A^{(1)}$. If $n>4$ and $v_2\in\C^n$ is the second basis from Section \ref{section CK for k-Dirac op}. Then with the usual convention $v_1\wedge v_2\in\Lambda^2\mF$ is a highest weight vector and $v_2.s=0$. This shows that $\mE\boxtimes\Lambda^2\mE\otimes\Lambda^2\mF\boxtimes\mF\boxtimes\Sp\subset A^{(1)}$. This proves the first part of the lemma. Now we use the Weyl dimension formula to compute the dimension of each module from the list.

The Weyl dimension formula works for semi-simple complex Lie algebras. We denote by $\Phi^+$ the set of all positive roots, by $\rho$ the lowest form. Let $\mV_\lambda$ be an irreducible module with highest weight $\lambda$. The Weyl dimension formula is
\begin{equation}\label{Weyl dimension formula}
 \dim(\mV_\lambda)=\frac{\prod_{\alpha\in\Phi^+}\langle\rho+\lambda,\alpha\rangle}{\prod_{\alpha\in\Phi^+}\langle\rho,\alpha\rangle}.
\end{equation}

We will use the same notation as in \cite{CSl}. Suppose that $n=2m$. Then $\Phi^+=\{e_i\pm e_j|1\le i<j\le m\}$ and $\rho=(m-1,m-2,\ldots,1,0)$. The denominator is 
$\prod_{\alpha\in\Phi^+}\langle\rho,\alpha\rangle=(2m-3)!(2m-5)!\ldots 3!1!(m-1)!.$
Suppose that $m\ge3$. The highest weight of $\mV_\lambda=S^3\mF\boxtimes\Sp_+$ is $(3+\frac{1}{2},\frac{1}{2},\ldots,\frac{1}{2})$. Thus the nominator is 
$\prod_{\alpha\in\Phi^+}\langle\lambda+\rho,\alpha\rangle=\frac{1}{6}(2m+1)!(2m-4)!(2m-6)!\ldots2!.$
So we get that $\dim(S^3\mF\boxtimes\Sp_+)=\frac{2^{m-2}}{3}(2m+1)(2m)(2m-1)$. The highest weight of $\mV_\mu=\mF\boxtimes\Lambda^2\mF\boxtimes\Sp_+$ is $\mu=(2+\frac{1}{2},1+\frac{1}{2},\frac{1}{2},\ldots,\frac{1}{2})$ Then
$ \prod_{\alpha\in\Phi^+}\langle\mu+\rho,\alpha\rangle=\frac{1}{6}(2m+1)(2m-1)!(2m-3)!(2m-6)!(2m-8)!\ldots2!.$
This gives that $\dim(\mV_\mu)=\frac{(2m+1)(2m-1)(2m-3)2^{m-1}}{3}$. 

Let us now consider the case $n=2m+1$. Then $\Phi^+=\{e_i\pm e_j,e_i|1\le i<j\le m\}$ and $\rho=(m-\half,m-\frac{3}{2},\ldots,\frac{3}{2},\half)$. The denominator is 
$ \prod_{\alpha\in\Phi^+}\langle\rho,\alpha\rangle=(2m-2)!(2m-4)!\ldots 2!(m-\half)(m-\frac{3}{2})\ldots\frac{3}{2}\half.$
Suppose that $m\ge2$. The highest weight of $\mV_\sigma=S^3\mF\boxtimes\Sp$ is $(3+\frac{1}{2},\frac{1}{2},\ldots,\frac{1}{2})$. Thus the nominator is 
$\prod_{\alpha\in\Phi^+}\langle\sigma+\rho,\alpha\rangle=\frac{1}{6}(2m+2)!(2m-3)!(2m-5)!\ldots3!1!.$
So we get that $\dim(S^3\mF\boxtimes\Sp)=\frac{2^m}{6}(2m+2)(2m+1)2m$. The highest weight of $\mV_\nu=\mF\boxtimes\Lambda^2\mF\boxtimes\Sp$ is $\nu=(2+\frac{1}{2},1+\frac{1}{2},\frac{1}{2},\ldots,\frac{1}{2})$. Then
$\prod_{\alpha\in\Phi^+}\langle\nu+\rho,\alpha\rangle=\frac{1}{3}(2m+2)(2m)!(2m-2)!(2m-5)!(2m-7)!\ldots3!1!.$
This gives that $\dim(\mV_\nu)=\frac{(m+1)m(m-1)2^{m+3}}{3}$.  $\Box$

Now we need to find the Cartan characters. 
\begin{lemma}\label{lemma Cartan char}
Let $k=2$ and $n\ge3$. Then the sequence of Cartan characters is equal to $(2n-2)s,(2n-3)s,\ldots,3s,2s,0,0,0$.
\end{lemma}
Proof: Let us first consider the case $n=3$. Let us choose the ordered basis $e_1\otimes\varepsilon_1,e_2\otimes\varepsilon_2,(e_1+e_2)\otimes\varepsilon_3,(e_1-e_2)\otimes\varepsilon_3,e_2\otimes\varepsilon_1,e_1\otimes\varepsilon_2$ of $\mV^\ast$. The corresponding affine coordinates are 
\begin{eqnarray}\label{affine coordinates1}
(t_1,\ldots,t_6)\mapsto
\left(
\begin{matrix}
t_1&t_5\\
t_6&t_2\\
t_3+t_4&t_3-t_4
\end{matrix}
\right).
\end{eqnarray}
We have to show that the corresponding Cartan characters are $8,6,4,0,0,0$. The proof of Lemma \ref{quadratic monogenic spinors lemma}. shows that the space of quadratic monogenic spinors is an irreducible $\lasl(2,\C)\times\laso(3,\C)$ isomorphic to $S^2\mE\otimes S^2_0\mF\boxtimes\Sp$. The complex dimension of this module is equal to $18$. From the description of this module given in Lemma \ref{quadratic monogenic spinors lemma}. follows that the first two Cartan characters are equal to $8,6$. It suffices to show that the last three Cartan characters are zero. That is: there is no monogenic quadratic spinor in the variables $t_4,t_5,t_6$.

The description of the basis of $S^2\mE\otimes S^2\mF$ given in Lemma \ref{quadratic monogenic spinors lemma}. shows that a polynomial $f\in S^2\mE\otimes S^2\mF\cap\C[t_4,t_5,t_6]$ is necessarily of the form $f=at_4^2+bt_5^2+ct_6^2$ for $a,b,c\in\C$. So we have to consider a spinor of the form $(a^\mu (e_1-e_2)\odot(e_1-e_2)\otimes\varepsilon_3\odot\varepsilon_3+b^\mu e_2\odot e_2\otimes\varepsilon_1\odot\varepsilon_1+c^\mu e_1\odot e_1\otimes\varepsilon_2\odot\varepsilon_2)\otimes s^\mu$ where $\mu=1,2$. This element belongs to the kernel of the map (\ref{description of prolongation}) iff all coefficients $a^\mu,b^\mu,c^\mu$ are zero. This proves the claim for $n=3$.

For $n>3$ we choose the ordered basis $e_1\otimes\varepsilon_1,e_2\otimes\varepsilon_1,e_1\otimes\varepsilon_2,\ldots,e_2\otimes\varepsilon_{n-3},e_1\otimes\varepsilon_{n-2},e_2\otimes\varepsilon_{n-1},(e_1+e_2)\otimes\varepsilon_n,(e_1-e_2)\otimes\varepsilon_n,e_2\otimes\varepsilon_{n-2},e_1\otimes\varepsilon_{n-1}$. The affine coordinates are then
\begin{eqnarray}\label{coordinates on affine space of matrices}
(t_1,t_2,t_3,\ldots,t_{2n})\mapsto
\left(
\begin{matrix}
t_1&t_2\\
t_3&\ldots\\
\ldots&\ldots\\
t_{2n-5}&t_{2n-1}\\
t_{2n}&t_{2n-4}\\
t_{2n-3}+t_{2n-2}&t_{2n-3}-t_{2n-2}
\end{matrix}
\right).
\end{eqnarray}
The claim follows from the description of the space of quadratic monogenic spinors in the proof of Lemma \ref{quadratic monogenic spinors lemma}. and the case $n=3$. $\Box$

%We will need later on the following equality.
%\begin{lemma}
%\sum_{i=1}^n i(n+1-i)={n+2\choose 3}
%\end{lemma}
%Proof: The left hand side 
%\begin{eqnarray*}
%\sum_{i=1}^n i(n+1-i)&=&(n+1)(1+2+\ldots+n)-\sum_{i=1}^ni^2\\
%&=&(n+1){n+1\choose2}-\frac{n(n+1)(2n+1)}{6}\\
%&=&\frac{n(n+1)}{6}(3n+3-2n-1)\\
%&=&{2n+2\choose3}.
%\end{eqnarray*}
%We have used here that $\sum_{i=1}^ni^2=\frac{n(n+1)(2n+1)}{6}$ which can be shown by induction. $\Box$

\begin{thm}\label{thm involutivity Dirac op}
The first prolongation of the tableau associated to the 2-Dirac operator is involutive.
\end{thm}
Proof: By the previous lemma the right hand side of the Cartan test (\ref{Cartan test1}) is equal to
%\begin{eqnarray*}
$s\sum_{i=1}^{2(n-1)}i(2n-1-i)-2s(n-1) %&=&s{2(n-1)+2\choose 3}-2s(n-1)\\
=s{2n\choose3}-2s(n-1).$
%\end{eqnarray*}
We used that $\sum_{i=1}^n i(n+1-i)={n+2\choose 3}$. Now we use the lower bound on $\dim(A^{(1)})$ from Lemma \ref{lemma second prolongation}. Recall that $s=\dim(\Sp_+\oplus\Sp_-)=2^m$ where $n=2m$ if $n$ is even while $n=2m+1$ if $n$ is odd. For $n\ge5$ we have that 
%\begin{eqnarray*}
 $\dim(A^{(1)})\ge \frac{2^{m+1}(n+1)n(n-1)+2^{m-1}(n+1)(n-1)(n-3)}{3}
 %2.4\frac{2^{m-2}}{3}(2m+1)(2m)(2m-1)+2.2\frac{(2m+1)(2m-1)(2m-3)2^{m-1}}{3}\\
 %&=&\frac{2^m}{6}((4m+2)(4m-2)(4m-3))\\
 %&=&\frac{2^m}{6}((4m-2)(16m^2-4m-6))\\
 %&=&\frac{2^m}{6}((4m-2)(4m-1)4m-6(4m-2))\\
 =2^m{2n\choose 3}-2.2^m(n-1).$
%\end{eqnarray*}
We have equality in the Cartan test and the tableau is involutive. Let us consider $n=4$. The module $S^3_0\mF\boxtimes\Sp_+$ is isomorphic to 
$S^3\C^2\otimes S^3\C^2\boxtimes\Sp_+\cong S^4\C^2\otimes S^3\C^2$. The dimension is equal to $20$. %Recall that $\laso(4)\cong\laso(3)\oplus\laso(3)$ and that the defininig representation 
%$\C^4$ is isomorphic to $\C^2\otimes\C^2$. First let us consider the piece $(2,1)\otimes(\C^2\otimes S^4\C^2\oplus S^4\C^2\otimes\C^2)$. We have $\dim(2,1)=2,\dim(S^4\C^2)=5$. % 
The dimension of the latter piece from Lemma \ref{lemma second prolongation} is clearly $40$. Since $\dim S^3\mE=4$ we get that $\dim(A^{(1)})\ge4.40+40=200$. On the other hand the sum of Cartan characters is equal to $4{8\choose3}-2.4.3=200$ and this is again an involutive tableau. This completes the proof for $n$ even.
%The case $n=2m+1$ for $m\ge2$ is 
%\begin{eqnarray*}
 %\dim(A^{(1)})&\ge&4\frac{2^m}{6}(2m+2)(2m+1)2m+2\frac{(m+1)m(m-1)2^{m+3}}{3}\\
 %&=&\frac{2^m}{6}(4(2m+2)2m(2m+1+2m-2))\\
  %&=&\frac{2^m}{6}((4m+4)4m(4m-1))\\
  %&=&\frac{2^m}{6}((16m^2+12m-4)4m)\\
  %&=&\frac{2^m}{6}((4m+2)(4m+1)4m-6(4m))\\
  %&=&2^m{2n\choose3}-2.2^m(n-1).
%\end{eqnarray*}
%Since $s=2^m$ the claim follows. 
The last remaining case is $n=3$. Recall that $\lasl(2,\C)\cong\laso(3,\C)$ and $\Sp\cong\C^2,\mF\cong\lasl(3,\C)$. Then $\dim(S^3\mE\otimes S^3_0\mF\boxtimes\Sp)=32$. The sum of the Cartan characters is $2(4+2.3+3.2)=32$.  $\Box$

%The last remaining case is $n=5$. The dimension of the module $S^3\mE\otimes S^3\mF\boxtimes\Sp$ is equal to $320$ as we have already computed above. Now $\Lambda^2\mF$ has highest weight $(1,1)$ and $\mF\boxtimes\Lambda^2\mF\boxtimes\Sp$ has highest weight $(\frac{5}{2},\frac{3}{2})$. The Weyl character formula shows that the dimension of this module is equal to $\frac{2.6.4.2}{2.\frac{3}{2}\frac{1}{2}}=64$. Thus $\dim(A^{(1)})=320+128=448$. Thus sum is $4{10\choose 3}-2.4.4=448$. This completes the proof of the theorem. 

\subsection{Initial conditions for the 2-Dirac operator}\label{section initial conditions}
We now use the Cartan-K\"ahler theorem to characterize the set of initial conditions for the 2-Dirac operator. Let us first consider the case $n=3$. Already on the lowest dimensional case we can illustrate the power of the Cartan-K\"ahler theorem. The general case will be given below.

Let us recall that we chose in the proof of Lemma \ref{lemma Cartan char}. an ordered basis of $M(3,2,\R)$ with affine coordinates in (\ref{affine coordinates1}).
%\begin{eqnarray}
%(t_1,\ldots,t_6)\mapsto
%\left(
%\begin{matrix}
%t_1&t_5\\
%t_6&t_2\\
%t_3+t_4&t_3-t_4
%\end{matrix}
%\right).
%\end{eqnarray}
Let us now denote the natural coordinates on the space $J^2\Sp$ of 2-jets of spinors by $\{t_i,s^\mu,u^\mu_j,p^\mu_{ij}\}$ so that the canonical contact system is $\theta^\mu=ds^\mu-u^\mu_j\omega_j,\theta^\mu_i=du^\mu_i-p^\mu_{ij}\omega_j$ where $\omega_i=dt_i,p^\mu_{ij}=p^\mu_{ji},\mu=1,2,i,j=1,\ldots,6$. Let $\Sigma$ be the subset of $J^2\Sp$ of 2-jets of monogenic spinors. Then local solutions of $\partial f=0$ are in 1-1 correspondence with integral manifols in $\Sigma$ of the pullback of the canonical system satisfying the independence condition given by $\omega_1,\ldots,\omega_6$. On such manifold we can consider $s^\mu,u^\mu_j,p^\mu_{ij}$ as functions of $t_1,\ldots,t_6$.

For $\mu=1,2,i=1,2,3$ and $j=i,i+1,\ldots,4$ let $f^\mu_{ij}$ be arbitrary real analytic functions of variables $t_1,\ldots,t_i$. Now the proof of the Cartan-K\"ahler theorem gives that there is a unique integral manifold of the canonical system satisfying the independence condition passing through the point $t_1=\ldots=t_6=s^\mu=u^\mu_1=\ldots=u^\mu_6=0$ such that the following set of equations
\begin{eqnarray*}\label{initial conditions}
 p^\mu_{1j}(t_1,0,0,0,0,0)&=&f^\mu_{1j}(t_1);j=1,2,3,4\\
 p^\mu_{2j}(t_1,t_2,0,0,0,0)&=&f^\mu_{2j}(t_1,t_2);j=2,3,4\\
 p^\mu_{3j}(t_1,t_2,t_3,0,0,0)&=&f^\mu_{3j}(t_1,t_2,t_3);j=3,4
\end{eqnarray*}
holds on the integral manifold. Recall that this gives existence of a monogenic spinor which satisfies the system of initial conditions (\ref{initial conditions}) on an open neighbouhood of the given point.

Since the $k$-Dirac operator is a constant coefficient system it suffices to understand homogeneous parts of monogenic spinors. The system of equations (\ref{initial conditions}) is equivalent to the following: given arbitrary homogeneous spinors $f_1,f_2$ of homogeneity $r$, resp. $r-1$ where $r\ge2$ in variables $t_1,t_2,t_3$ then there is a unique monogenic spinor $f_1+t_4f_2+g$ where $g$ is a homogeneous spinor of degree $r$ such that the sum of degrees of the variables $t_4,t_5,t_6$ of each monomial appearing in a component of $g$ is at least equal to 2.

For example consider quadratic monogenic spinors. The space of quadratic monogenic spinors is naturally isomorphic to the vector space $\{a^\mu_{ij}t_it_j+b^\mu_lt_4t_l\}$ where $a^\mu_{ij}\in\C,b_l^\mu\in\C$ are arbitrary constants symmetric in $i,j,l=1,2,3$. Note that the dimension of this space is $2({4\choose2}+3)=18$ which agrees with the previous computations. Cubic monogenic spinors are naturally isomorphic to the space $\{a^\mu_{ijl}t_it_jt_l+b^\mu_{uv}t_4t_ut_v\}$ where $a^\mu_{ijl}\in\C,b_{uv}^\mu\in\C$ are arbitrary constants symmetric in $i,j,l,u,v=1,2,3$.  The dimension of the space of these coefficients is $2({5\choose 3}+{4\choose 2})=32$. 

For general $n$ we use the same coordinates as in (\ref{coordinates on affine space of matrices}). Let $g_1,g_2$ be homogeneous spinors on $M(n,2,\R)$ of homogeneity $r$, resp. $r-1$ with $r\ge2$ in variables $t_1,\ldots,t_{2n-3}$. Then there is a unique monogenic spinor $g_1+t_{2n-2}g_2+g$ where $g$ is a homogeneous spinor of degree $r$ such that the sum of degrees of the variables $t_{2n-2},t_{2n-1},t_{2n}$ of each monomial appearing in a component of $g$ is at least equal to 2. We have proved the following theorem.

\begin{thm}\label{thm initial conditions}
The vector space of homogeneous monogenic spinors of degree $r\ge2$ for the 2-Dirac operator in dimension $n\ge3$ is naturally isomorphic to the direct sum of vector spaces of homogeneous spinors of degree $r$ and $r-1$ in the variables $t_1,t_2,\ldots,t_{2n-3}$ from (\ref{coordinates on affine space of matrices}).
\end{thm}

\section{The Cartan-K\"ahler theorem for the parabolic $k$-Dirac operator}

\subsection{The canonical linear Pfaffian system on $J^2\Sp$}
Let us recall that we are working on the affine set $\osu$ from Section \ref{section parabolic D operator}. We write coordinates on $M(n,k,\R)$ as $x_{\alpha i}$ and on $A(k,\R)$ as $y_{rs}$. Then $\partial_{\alpha i}$ and $\partial_{rs}=-\partial_{sr}$ stand for the coordinate vector fields. We use convention $dy_{rs}(\partial_{ij})=\delta_{ir}\delta_{js}-\delta_{jr}\delta_{is}$. Then we can write $dy_{rs}=-dy_{sr}$. We set $L_{\alpha i}=\partial_{\alpha i}-\frac{1}{2}x_{\alpha j}\partial_{ij}$ and call them left invariant vector fields.

As in the case of the Euclidean $2$-Dirac operator we will need the first prolongation of the canonical linear Pfaffian system living on the space of 1-jets of monogenis spinors ,i.e. we will need to work on the space of 2-jets of monogenic spinors on $\osu$. We write coordinates on the space $J^2\Sp$ of 2-jets of spinors over $\osu$ as $\{x_{\alpha i},y_{rs},s^\mu,u^\mu_{\alpha i},v^\mu_{rs},a^\mu_{\alpha i\beta j},b^\mu_{\alpha i rs},c^\mu_{rs uv}\}$ with the relations $v^\mu_{rs}=v^\mu_{sr},a^\mu_{\alpha i\beta j}=a^\mu_{\beta j\alpha i},b^\mu_{\alpha irs}=-b^\mu_{\alpha isr},c^\mu_{rs uv}=c^\mu_{uv rs}=-c^\mu_{sr uv}=-c^\mu_{rs vu}$. The canonical Pfaffian system $\cI$ is generated by 1-forms
$\theta^\mu=ds^\mu-u^\mu_{\alpha i}dx_{\alpha i}-\frac{1}{2}v^\mu_{rs}dy_{rs}\label{canonical LPS I},
\theta^\mu_{\alpha i}=du_{\alpha i}^\mu-a^\mu_{\alpha i\beta j}dx_{\beta j}-\frac{1}{2}b^\mu_{\alpha irs}dy_{rs},
\theta^\mu_{uv}=dv_{uv}^\mu-b^\mu_{\alpha i uv}dx_{\alpha i}-\frac{1}{2}c^\mu_{uvrs}dy_{rs}.$
Here we are summing over all $r,s=1,\ldots,k$ and so the factor $\frac{1}{2}$ appears there. We denote the ideal generated by these forms by $I$. 

We will need to introduce new coordinates which are more adapted for the operator $D$. In the first place we have to find the dual 1-forms to the vector fields $L_{\alpha i},\partial_{rs}$, i.e. we are looking for 1-forms such that
$\omega_{\alpha i}(L_{\beta j})=\delta_{\alpha\beta}\delta_{ij},\omega_{\alpha i}(\partial_{rs})=\omega_{rs}(L_{\alpha i})=0,\omega_{rs}(\partial_{ij})=\delta_{ir}\delta_{js}-\delta_{is}\delta_{jr}.$ These forms will give for each $x\in\osu$ an isomorpism $T^\ast_x\osu\cong\mE\otimes\mF\oplus\Lambda^2\mE$ where $\Lambda^2\mE$ is isomorphic to the span of all $(\omega_{rs})_x$ and $\mE\otimes\mF$ is the span of all $(\omega_{\alpha i})_x$. We will not carefully distinguish between complex and real representations as we did in the previous sections. The meaning should be clear from the context. We find that $\omega_{\alpha i}=dx_{\alpha i}$ and $\omega_{rs}=dy_{rs}-\frac{1}{2}(x_{\beta r}dx_{\beta s}-x_{\beta s}dx_{\beta r})$. We have $d\omega_{\alpha i}=0,d\omega_{rs}=\sum_\alpha\omega_{\alpha s}\wedge\omega_{\alpha r}$. 

Substituting $\omega_{\alpha i},\omega_{rs}$ into the formula for $\theta^\mu$ we obtain
$\theta^\mu=ds^\mu-\sigma^\mu_{\alpha i}\omega_{\alpha i}-\frac{1}{2}v^\mu_{rs}\omega_{rs}$
where $\sigma^\mu_{\alpha i}=u^\mu_{\alpha i}-\frac{1}{2}x_{\alpha j}v^\mu_{ij}$. We set
$A^\mu_{\alpha i\beta j}=a^\mu_{\alpha i\beta j}-\frac{1}{2}(x_{\alpha s} b^\mu_{\beta jis}+x_{\beta t}b_{\alpha ijt}^\mu)+\frac{1}{4}x_{\alpha s}x_{\beta t}c^\mu_{ri sj}-\frac{1}{2}\delta_{\alpha\beta}v^\mu_{ij},
B^\mu_{\alpha i js}=b_{\alpha ijs}^\mu-\frac{1}{2}x_{\alpha t}c^\mu_{itjs},C^\mu_{rs kl}=c^\mu_{rskl}.$
Then $A^\mu_{\alpha i\beta j}-A^\mu_{\beta j\alpha i}=\delta_{\alpha\beta}v^\mu_{ij}$. This is compatible with (\ref{Lie bracket of left invariant vector fields}). The forms $\theta^\mu_{\alpha i},\theta^\mu_{rs}$ are then
$\theta^\mu_{\alpha i}=d\sigma^\mu_{\alpha i}+\frac{1}{2}x_{\alpha j}\theta_{ij}-A^\mu_{\alpha i\beta j}\omega_{\beta j}- B^\mu_{\alpha irs}\omega_{rs},
\theta^\mu_{rs}=dv^\mu_{rs}-B^\mu_{\beta j rs}\omega_{\beta j}-\frac{1}{2}C^\mu_{rsuv}\omega_{uv}.$

\subsection{Vanishing of torsion}\label{section killing torsion}
In this section we argue that the torsion of the linear Pfaffian system associated to the $k$-Dirac operator and to its prolongation vanishes. We state a necessary lemma from \cite{S} and set the notation. We define a grading on the space of polynomials $\C[x_{\alpha i},y_{rs}]$. The \textit{weighted degree} of linear polynomials is $deg_w(x_{\alpha i})=1,deg_w(y_{rs})=2$. We extend this to the space of monomials such that $deg_w$ is a morphism of $(\C[x_{\alpha i},y_{rs}],.)\ra(\Z,+)$. Then $deg_w(f)=r$ iff $f$ is a sum of monomials of weighted degree $r$. We say that a spinor $\psi$ on $\osu$ is of weighted degree $r$ if each component of $\psi$ is a weighted polynomial of degree $r$ in the preferred trivialization.
\begin{lemma}\label{killing torsion}
 Let $\psi$ be a homogeneous monogenic spinor (in the Euclidean setting) of degree $r$ on $M(n,k,\R)$, i.e. $\partial\psi=0$. Let $g\in\C[y_{rs}]$ be an arbitrary homogeneous polynomial of degree $l$. Then there is a parabolic monogenic spinor $\Psi$ homogeneous of weighted degree $r+2l$ on $\osu$, i.e. $D\Psi=0$, which is of the form $\Psi=g\psi+l.o.t.$ where $l.o.t.$ stands for a spinor on $\osu$ whose components are polynomials which are of degree strictly smaller than $l$ in $y$-variables.
\end{lemma}
Proof: This is Lemma 8.6.2 from \cite{S}. $\Box$

Let us choose $k\in\{1,2\}$. We denote the space of $k$-jets of monogenic spinors by $\Sigma$. The lemma implies that there is an integral manifold passing through any point in the fibre of the canonical projection $\Sigma\ra\osu$ over the origin $0\in\osu$. The tangent space of the integral manifold is an integral element and so by Lemma \ref{lemma vanishing of torsion}. the torsion $[T]=0$ vanishes identically in the fibre over $0\in\osu$. 

The flow of (the projection of) a right invariant vector field $X$ on the homogeneous space is symmetry of the operator $D$. By flows of such vector fields we can move any point $x\in\osu$ to any given point $x'\in\osu$. The induced action on $\Sigma$ is compatible with the canonical projection to $\osu$. The flow of the field $X$ preserve the ideal $I$ and thus also the tableau and the torsion is invariant along the flow lines. Since the torsion vanishes in the fibre over $0\in\osu$ it has to vanish everywhere on $\Sigma$.

\subsection{Non-involutivity of the tableau associated to $k$-Dirac operator $D$}

The space of 1-jets of spinors on $\osu$ is the set $J^1\Sp=\{(x_{\alpha i},y_{rs},s^\mu,\sigma^\mu_{\alpha i},v^\mu_{rs})\}$ with canonical linear Pfaffian system generated by the forms $\theta^\mu$ and the indepence condition $\omega_{\alpha i},\omega_{rs}$. The structure equations are 
$d\theta^\mu=-d\sigma^\mu_{\alpha i}\wedge\omega_{\alpha i}-\frac{1}{2}dv^\mu_{rs}\wedge\omega_{rs}-\frac{1}{2}v^\mu_{rs}d\omega_{rs}.$
We use abstract index notation and the Einstein summation convention. We can then write $(\varepsilon_\alpha.):\Sp\ra\Sp,(\varepsilon.s)^\mu=(\varepsilon.)^\mu_\nu s^\nu$ for any spinor $s^\nu\in\Sp$ and $\varepsilon\in\mF$.
 
Then a 1-jet from $J^1\Sp$ is a 1-jet of monogenic spinor iff $\sum_\alpha(\varepsilon_\alpha)^\nu_\mu\sigma^\mu_{\alpha i}=0$ for all $i=1,\ldots,k$. We may take $\pi^\varepsilon$ to be the forms $d\sigma^\mu_{\alpha i}$ with $\alpha>1$ and $dv^\mu_{rs}$ with $r<s$. For each $x\in\osu:\mV^\ast_x\cong\mE\otimes\mF\oplus\Lambda^2\mE$ and $\mW_x\cong\Sp$. The tableau is at any point isomorphic to $\mE\otimes\T\oplus\Lambda^2\mE\otimes\Sp$ while the torsion is represented by $[-\frac{1}{2}v^\mu_{rs}d\omega_{rs}]$. From the discussion in the previous section follows that the torsion vanishes identically.

The Cartan characters with respect to the ordered basis $e_1\otimes\varepsilon_1,\ldots,e_1\otimes\varepsilon_k,\ldots,e_{n-1}\otimes\varepsilon_1,e_{n-1}\otimes\varepsilon_k,e_1\wedge e_2,\ldots,e_{k-1}\wedge e_k,e_n\otimes\varepsilon_1,\ldots,e_n\otimes\varepsilon_k$ of $\mV^\ast$ are $s_i=s,s_j=0$ for $i\le k(n-1)+{k\choose2}<j$ and so the right hand side of (\ref{Cartan test1}) is equal to
$ %s_1+2s_2+\ldots+(nk+{k\choose2})s_{nk+{k\choose2}}%&=&s(1+2+\ldots+(k(n-1)+{k\choose2}))\nonumber\\
s{k(n-1)+{k\choose2}+1\choose2}.$
The first prolongation is clearly isomorphic to 
\begin{equation}\label{first prolongation of parabolic D operato}
A^{(1)}\cong \mM\oplus\Lambda^2\mE\otimes\mE\otimes\T\oplus S^2(\Lambda^2\mE)\otimes\Sp
\end{equation}
where $\mM$ is the space of the quadratic monogenic spinors (in the Euclidean setting) described in the proof of Lemma \ref{quadratic monogenic spinors lemma}.  The dimension of the prolongation is $\dim(A^{(1)}) %=s{k(n-1)+1\choose 2}-s{k\choose2}+s{k\choose 2}k(n-1)+s{{k\choose2}+1\choose2}
=s{k(n-1)+{k\choose2}+1\choose2}-s{k\choose2}.$

We see that we do not have equality in the Cartan test (\ref{Cartan test1}) and thus the tableau is not involutive. We have to prolong this system as we in the case of the 2-Dirac operator $\partial$. The interpretation of the tableau and its prolongations is the following. Let $J^i_x\cM$ be the space of $i$-jets of monogenic spinors at a point $x\in\Sigma$. The tableau is isomorphic to the kernel of the canonical projection $J^1_x\cM\ra J^0_x\cM$, the first prolongation of the tableau is isomorphic to the kernel of  $J^2_x\cM\ra J^1_x\cM$, the prolongation of the first prolongation is then isomorphic to the kernel of $J^3_x\cM\ra J^2_x\cM$ and so on. 

\subsection{Involutivity of the first prolongation of the parabolic 2-Dirac operator $D$}
The structure equations on $J^2\Sp$ are 
$d\theta^\mu=0,
d\theta^\mu_{\alpha i}=\frac{1}{2}x_{\alpha j}d\theta_{ij}-dA^\mu_{\alpha i\beta j}\wedge\omega_{\beta j}-\frac{1}{2}dB^\mu_{\alpha irs}\wedge\omega_{rs}-\frac{1}{2}B^\mu_{\alpha irs}d\omega_{rs},
d\theta^\mu_{rs}=-dB^\mu_{\beta j rs}\wedge\omega_{\beta j}-\frac{1}{2}dC^\mu_{rs uv}\wedge\omega_{uv}-\frac{1}{2}C^\mu_{rs uv}d\omega_{uv}$ all modulo $I$.

The space of 2-jets of monogenic spinors is the subset of $J^2\Sp$ where the following set of relations holds. In the first place: $\sum_\alpha(\varepsilon_\alpha.)^\mu_\nu\sigma^\nu_{\alpha i}=0,\sum_\alpha(\varepsilon_\alpha.)^\mu_\nu B^\nu_{\alpha irs}=0$ holds for all $i,r,s$. From the variables $A^\mu_{\alpha i\beta j}$ only those with $\alpha,\beta>1$ are free. There is one more system of equations $\sum_{\alpha,\beta>1}[(\varepsilon_\alpha.)^\mu_\rho,(\varepsilon_\beta.)^\rho_\nu]A^\nu_{\alpha i\beta j}=(-n+2)v^\mu_{ji}$. 

We find that for any $x\in\osu:\mV_\ast\cong\mE\otimes\mF,\mW\cong \mE\otimes\T\oplus\Lambda^2\mE\otimes\Sp$ and that the tableau $A$ is isomorphic to the first prolongation (\ref{first prolongation of parabolic D operato}) from the previous section. %$S^2\mE\otimes S^2\mF\boxtimes\Sp\oplus\Lambda^2\mE\otimes\Lambda^2\mF\boxtimes\Sp\oplus\Lambda^2\mE\otimes\mE\otimes\T\oplus S^2(\Lambda^2\mE)\otimes\Sp$. 
The torsion vanishes identically on $\Sigma$ by the same argument as in the previous section.

\begin{lemma}
The Cartan characters of the tableau $A$ are $(2n-1)s,(2n-2)s,(2n-3)s,\ldots,3s,2s,0,0,0$.
\end{lemma}
Proof: Let us choose the origin $x=0\in\Sigma$. Let us order basis of $\mV^\ast$ by putting the vector $e_1\wedge e_2\in\Lambda^2\mE$ in the first place and then we put the basis of $\mE\otimes\mF$ ordered in the same way as in Lemma \ref{lemma Cartan char}. Then the first Cartan character is equal to the dimension of $S^2(\Lambda^2\mE)\otimes\Sp\oplus\Lambda^2\mE\otimes\mE\otimes\T$. This number is equal to $s(1+2(n-1))$. The other Cartan characters clearly coincide with the Cartan characters from Lemma \ref{lemma Cartan char}. $\Box$

\begin{lemma}
The first prolongation $A^{(1)}$ is isomorphic to the direct sum of the corresponding irreducible $\lasl(2,\C)\oplus\laso(n)$-modules from the table (\ref{ternary monogenic spinors}) and $\Lambda^2\mE\otimes\mM\oplus S^2(\Lambda^2\mE)\otimes\mE\otimes\T\oplus S^3(\Lambda^2\mE)\otimes\Sp$ where $\mM$ is the module isomorphic the space of quadratic monogenic spinors (in the Euclidean setting) from Lemma \ref{quadratic monogenic spinors lemma}.
\end{lemma}
Proof: Follows by the definition of $A^{(1)}:=\mV^\ast\otimes A\cap S^2\mV^\ast\otimes\mW$. $\Box$

\begin{thm}\label{thm involutivity D op}
The tableau of the first prolongation of the parabolic 2-Dirac operator is involutive.
\end{thm}
Proof: The right hand side of the Cartan test (\ref{Cartan test1}) is
$ %s_1+2s_2+\ldots+2s_(2n-2)=
s\sum_{i=1}^{2n-1}i(2n-i)-(2n-1)s
%&=&s{2n+1\choose3}-s(2n-1)\\
=s\frac{(2n-1)}{6}(4n^2+2n-6).$
The left hand side is equal to
$\dim(A^{(1)})=s{2n\choose3}-2s(n-1)+s({2(n-1)+1\choose 2}-1)+2s(n-1)+s
%&=&s(2n-1)\frac{4n^2-4n+6n-6}{6}\\
=s(2n-1)\frac{4n^2+2n-6}{6}.$
Here we have used that $\dim(\mM)=s({2(n-1)+1\choose 2}-1)$ proved in Lemma \ref{thm involutivity Dirac op}. and that the dimension of the space of cubic monogenic spinors (in the Euclidean setting) is $s{2n\choose3}-2s(n-1)$ which was shown in the proof of Theorem \ref{thm involutivity Dirac op}.  $\Box$

\end{document}